\documentclass[12pt]{amsart}

\usepackage[latin1]{inputenc}
\usepackage{amsfonts,amsmath,amssymb,amsthm}

\newcommand{\Nn}{\mathbb{N}}

\newcommand{\Qq}{\mathbb{Q}} 

\renewcommand{\epsilon}{\varepsilon}
\renewcommand{\le}{\leqslant}
\renewcommand{\ge}{\geqslant}
\renewcommand{\leq}{\leqslant}
\renewcommand{\geq}{\geqslant}

\newcommand{\defi}[1]{\emph{#1}}
\newcommand{\xbar}{\underline{x}}
\newcommand{\lambdabar}{{\underline{\lambda}}}

\newcommand{\Klambdabar}{{\overline{K(\underline{\lambda})}}}
\newcommand{\Klambdar}{{\overline{K({\lambda})}}}
\newcommand{\charac}{\mathop{\mathrm{char}}\nolimits}
\newcommand{\point}{^\ast}
\newcommand{\reduced}{'}

{\theoremstyle{plain}
\newtheorem{theorem}{Theorem}[section]    
\newtheorem*{addendumA}{Theorem \ref{th:mon}}    
\newtheorem*{addendumB}{Theorem \ref{th:mon}}   

\newtheorem{lemma}[theorem]{Lemma}       
\newtheorem{proposition}[theorem]{Proposition}  
\newtheorem{corollary}[theorem]{Corollary}   
}
{\theoremstyle{remark}
\newtheorem{definition}[theorem]{Definition}      
\newtheorem{remark}[theorem]{Remark}   
\newtheorem{example}[theorem]{Example}
}

\begin{document}

\title{Irreducibility of hypersurfaces}

\author{Arnaud Bodin}
\email{Arnaud.Bodin@math.univ-lille1.fr}

\author{Pierre D\`ebes}
\email{Pierre.Debes@math.univ-lille1.fr}

\author{Salah Najib}
\email{Salah.Najib@math.univ-lille1.fr}

\address{Laboratoire Paul Painlev\'e, Math\'ematiques, Universit\'e 
Lille 1, 59655 Villeneuve d'Ascq Cedex, France}



\date{\today}
 
\begin{abstract}

  Given a polynomial $P$ in several variables over an 
  algebraically closed
  field, we show that except in some special cases that we fully
  describe,  if one coefficient is allowed to vary, then the polynomial
 is irreducible for all but at most $\deg(P)^2-1$ values of the coefficient.
 We more generally handle the situation where several specified
 coefficients vary. 
\end{abstract}


\maketitle

\section{Introduction}

Classically  polynomials in $n\geq 2$ variables are 
generically absolutely irreducible: if the coefficients, in some algebraically 
closed ground field $K$, are moved a little bit but stay away from  
some proper Zariski closed subset, 
then the resulting polynomial is irreducible over $K$.
This is no longer true if only one specified coefficient is allowed
to vary. For example 
however one moves a non-zero coefficient of some homogeneous 
polynomial $P(x,y)\in K[x,y]$ of degree $d\geq 2$, it remains reducible 
over $K$. Yet it seems that this case is exceptional and 
that most polynomials are irreducible up
to moving any fixed coefficient away from finitely many values.  This
paper is aimed at making this more precise.

\subsection{The problem}
The problem can be posed in general as follows: given an algebraically
closed field $K$ (of any characteristic) and a polynomial $P\in K[\xbar]$
(with $\xbar= (x_1,\ldots,x_n)$), describe the ``exceptional'' {\textit{reducibility 
monomial sites} of $P$, that is those sets $\{Q_1,\ldots,Q_\ell \}$ of
monomials in $K[\xbar]$ for which $P+\lambda_ 1Q_1 +  \cdots +
\lambda_\ell  Q_\ell$ is \textit{generically reducible}, \textit{i.e.} reducible
in $\Klambdabar[\xbar]$\footnote{Given a field $k$, we denote by $\overline k$
an algebraic closure of $k$.},  where $\underline{\lambda} = 
(\lambda_1,\ldots,\lambda_\ell)$ is a $\ell$-tuple of independent 
indeterminates. When this is not the case, it follows from the Bertini-Noether theorem
that the polynomial with shifted coefficients 
$P+\lambda\point_ 1Q_1 + \cdots +\lambda\point_\ell 
Q_\ell$ is irreducible in $K[\xbar]$ for all $\lambdabar\point =
(\lambda\point_1,\ldots,\lambda\point_\ell)$ in a non-empty Zariski 
open subset of $K^\ell$ (and the converse is true).

The situation $\ell=1$ has been extensively studied in the literature, notably
for $Q_1=1$, that is when it is the constant term that is moved: see 
works of Ruppert \cite{Ru}, Stein \cite{St}, Ploski \cite{Pl}, Cygan \cite{Cy},
Lorenzini \cite{Lo}, Vistoli \cite{Vi}, Najib \cite{Na}, Bodin
\cite{Bo} et al. The central result in this case, which is known as
Stein's theorem, is that $P+\lambda$ is generically irreducible if and only
if $P(\xbar)$ is not a composed polynomial\footnote{
that is, is not of the form $r(S(\xbar))$ with $S\in K[\xbar]$ and $r\in K[t]$
with $\deg(r) \geq 2$.} (some say ``indecomposable'');
furthermore, the so-called spectrum of $P$ 
consisting of all $\lambda\point\in K$ such
that $P+\lambda\point$ is reducible in $K[\xbar]$, which from Bertini-Noether 
is finite in this case, is of cardinality $< \deg(P)$. This was first established by Stein in two
variables and in characteristic $0$, then extended to all characteristics by 
Lorenzini and finally generalized to $n$ variables by Najib. 
The result also extends to arbitrary monomials $Q_1$, and in fact
to arbitrary polynomials \cite{Lo} \cite{Bo}; the indecomposability assumption 
should be replaced  by the condition that $P/Q_1$ is not a composed rational 
function, and the bound $\deg(P)$ by $\deg(P)^2$.


\subsection{Our results}
We fully describe the reducibility monomial sites of polynomials in
the general situation $\ell \geq 1$ (theorem \ref{th:mon}). 
We obtain simple criteria for generic irreducibility, more practical
than the previous indecomposability type conditions. These results 
can be combined with some $\ell$-dimensional Stein-like
description of the irreducibility set (proposition \ref{prop:specialisation}).
Our contribution can be illustrated by the following three consequences. 

Recall $K$ is an algebraically closed field of any characteristic. Below
by Newton  representation of a polynomial in $n$ variables we merely mean 
the subset of all points $(a_1,\ldots,a_n)\in \Nn^n$ such that the monomial 
$x_1^{a_1} \cdots x_n^{a_n}$ appears in the polynomial with a non-zero 
coefficient.

\begin{theorem} \label{thm:typical}
  Let $P(\xbar) \in K[\xbar]$ be a non constant polynomial and
  $Q(\xbar)$
  be a monomial of degree $\leq
  \deg(P)$ and relatively prime to $P$. Assume that the monomials of
  $P$ together with $Q$ do not lie on a line in their Newton 
  representation\footnote{The result also holds if  $P$ is a monomial (in which
  case $P$ and $Q$ are lined up in the Newton representation).}
  and that 
    $Q$ is not a pure
  power\footnote{We say a polynomial $R\in K[\xbar]$ is a \defi{pure power} if 
  there exist
    $S\in K[\xbar]$ and $e>1$ such that $R = S^e$. The monomial
    $Q(\xbar)=x_1^{e_1}\cdots x_n^{e_n}$ is not a pure power if and
    only if $e_1, \ldots, e_n$ are relatively prime. }  in $K[\xbar]$.  Then 
    $P+\lambda Q$ is generically irreducible and the set
  of all $\lambda\point \in K$ such that $P+\lambda\point Q$ is
  reducible in $K[\xbar]$ is finite of cardinality $< 
  \deg(P)^2$. 
\end{theorem}

In particular a polynomial can always be made irreducible by changing only
one of its coefficients provided it is not divisible by a non-constant monomial.

The assumption on the monomials of $P$ and $Q$ is here to avoid
what we call the exceptional homogeneous case, that is, that  $P$ be of the form 
$h(m_1,m_2)$ with $h\in K[u,v]$ homogeneous and 
$m_1$, $m_2$ two monomials of degree $<\deg(P)$,
in which case for any monomial $Q=m_1^k m_2^{d-k}$ ($0\leq k\leq d=\deg(h)$), 
$P+\lambda Q$ is generically reducible. 

Pure power monomials $Q$, \textit{e.g.} $Q=1$, should also be excluded in 
theorem \ref{thm:typical}, but can nevertherless be dealt with under a slightly 
more general condition.

\begin{theorem} \label{thm:typical3}
  Let $P(\xbar) \in K[\xbar]$ be a non constant polynomial and
  $Q(\xbar)$ be a monomial of degree $\leq \deg(P)$ and relatively
  prime to $P$. Assume $P$ is not of the form
  $h(m,\psi)$ with $h\in K[u,v]$ an homogeneous polynomial, $m$ a
  monomial dividing $Q$ and $\psi\in K[\xbar]$ such that $\deg(P) >
  \max(\deg(m),\deg(\psi))$. Then 
    $P+\lambda Q$ is generically irreducible and the set of all 
    $\lambda\point \in K$
  such that $P+\lambda\point Q$ is reducible in $K[\xbar]$ is 
  finite
  and of cardinality $< \deg(P)^2$.
\end{theorem}

If $P$ is of the excluded form then, for $Q=m^{\deg(h)}$, the polynomial
$P+ \lambda Q$ is
generically reducible.

In the special case $Q=1$, the assumption on $P$ is that it is not of
the form $h(1,\psi)$ with $h\in K[u,v]$ homogeneous, $\deg_v(h) \geq 2$
and $\psi\in K[\xbar]$: this corresponds to the classical hypothesis
that $P$ is not a composed polynomial. Thus theorem \ref{thm:typical3}
is a generalization of Stein's theorem (except  for the bound which can 
be taken to be $\deg(P)$ in this special case).

As another typical consequence of our approach, we obtain that for 
$\ell \geq 2$, reducibility monomials are even more rare.

\begin{theorem} \label{thm:typical2}
  Let $P \in K[\xbar]$ be a non constant polynomial and, for $\ell \geq 2$, 
  $Q_1,\ldots,Q_\ell$ be $\ell$ monomials of degree
  $\leq \deg(P)$ and such that $P,Q_1,\ldots,Q_\ell$ are relatively
  prime.  Assume the monomials of $P$ together with $Q_1, \ldots, Q_\ell$ 
  do not lie on a line in their Newton representation. 
   If $\charac(K)=p>0$ assume further that at least one of 
  $P,Q_1,\ldots,Q_\ell$ is not a $p$-th power.  Then $P+\lambda_1 
  Q_1+\cdots + \lambda_\ell Q_\ell$ is generically irreducible and so
  $P+\lambda\point_1 Q_1+\cdots + \lambda\point_\ell Q_\ell$ is
  irreducible in $K[\xbar]$ for all $(\lambda\point_1,\ldots,\lambda\point_\ell)$ 
  in a non-empty Zariski open
  subset of $K^\ell$.\footnote{Prop. \ref{prop:specialisation} gives a more 
  explicit Stein-like description of the irreducibility set.} 
\end{theorem}

For example
$P(x_1,\ldots,x_n)+\lambda_1 x_1+\cdots + \lambda_n
x_n$ ($n\ge 2$) is generically irreducible.  
See corollary \ref{cor:funny} for further related results.


\subsection{Organization of the paper}
A starting ingredient of our method is the Bertini-Krull
theorem, which gives an iff condition for some polynomial $P+\lambda_1
Q_1+\cdots + \lambda_\ell Q_\ell$ to be generically irreducible. 
The Bertini-Krull theorem
is recalled in the preliminary section \ref{sec:uniqueness} which also
introduces some basic definitions used in the rest of the paper. We
also seize the opportunity to prove a useful uniqueness result
(theorem \ref{thm:unicity}) in the Bertini-Krull theorem, which to our
knowledge, was only known in the context of Stein's theorem.

Section \ref{sec:mon.sites} is the core of the paper. We investigate
the Bertini-Krull conclusion in the specific context of our problem to
finally obtain a general description of 
the reducibility monomial  sites of a given polynomial
(theorem \ref{th:mon}). Giving an exact description requires controlling
the possible overlaps of the special cases where reducibility monomial  sites
can exist. 
This comes down to proving (as in lemma \ref{lem:monunicity}) some uniqueness 
statements for  ``homogeneous decompositions'' of polynomials related to those studied in 
section  \ref{sec:uniqueness}.

Section \ref {sect:specialization} is devoted to
specializing the variables $\lambda_1,\ldots,\lambda_\ell$. 
For $\ell = 1$,  we use the generalization of Stein's theorem due 
to Lorenzini \cite{Lo} and Bodin \cite{Bo} to give an upper bound
for the cardinality of the set of  exceptional values $\lambda\point$ 
making $P+\lambda\point Q$ reducible in $K[\xbar]$.
A version of this estimate can be derived inductively 
for the situation $\ell \geq 1$, for which the classical Bertini-Noether theorem 
can also be used.  We then complete the proof of the results from the
introduction and give some further corollaries.


\subsection{Main Data and Notation}\label{subs:notation}

The following is given and will be retained throughout the paper:
\begin{itemize}
\item an algebraically closed field $K$ of characteristic $0$ or
  $p>0$,
  
\item an integer $\ell \geq 0$ and an $\ell$-tuple
  $\lambdabar=(\lambda_1,\ldots,\lambda_\ell)$ of independent
  variables (algebraically independent over $K$); for $\ell=0$, the
  convention is that no variable is given,
  
\item an integer $n\geq 2$ and an $n$-tuple $\xbar = (x_1,\ldots,x_n)$
  of new independent variables (algebraically independent over
  $\Klambdabar$),
  
\item $\ell +1$ distinct (up to multiplicative constants) non-zero
  polynomi\-als $P,Q_1,\ldots,Q_\ell \in K[\xbar]$ with
  $\max(\deg(P),\ldots,\deg(Q_\ell))>0$ and assumed further to be
  relatively prime if $\ell \geq 1$,
  
\item $F(\xbar,\lambdabar)= P(\xbar) + \lambda_1 Q_1(\xbar) + \cdots
  +\lambda_\ell Q_\ell(\xbar)$, which is an irreducible polynomial in
  $K[\xbar,\lambdabar]$ if $\ell \geq 1$.  (For $\ell \geq
  1$, $F(\xbar,\lambdabar)$ can be al\-ter\-na\-ti\-ve\-ly defined as a linear
  form in $(\lambda_0,\ldots,\lambda_\ell)$ (with $\lambda_0=1)$ with
  distinct non-zero and relatively prime coefficients in $K[\xbar]$).
\end{itemize}


\section{Around the Bertini-Krull theorem}\label{sec:uniqueness}

\subsection{Bertini-Krull theorem and homogeneous decompositions}

We start by recalling the Bertini-Krull theorem. We refer to 
  \cite[theorem 37]{Sc} where equivalence between conditions (1) and (4) 
  below is proved; equivalence between conditions (1), (2) and
  (3) is a special case of the standard Bertini-Noether theorem
  \cite[proposition 8.8]{FrJa}.
  

\begin{theorem}[Bertini, Krull]
\label{th:BK} In addition to \S \ref{subs:notation}, 
assume $\ell \geq 1$.  Then the following conditions are equivalent:
\begin{enumerate}
\item $F(\xbar,\lambdabar\point)$ is reducible in $K[\xbar]$ for all
  $\lambdabar\point \in K^\ell$ such that \hfill  \break
  $\deg(F(\xbar,\lambdabar\point)) = \deg_{\xbar}(F)$.
  
\item The set of $\lambdabar\point \in K^\ell$ such that
  $F(\xbar,\lambdabar\point)$ is reducible in $K[\xbar]$ is
  Zariski-dense.
  
\item $F(\xbar,\lambdabar)$ is reducible in $\Klambdabar[\xbar]$.

\item
  \begin{enumerate}
  \item \label{iit:charp} either $\charac K = p >0$ and
    $F(\xbar,\lambdabar) \in K[\xbar^p,\lambdabar]$, where $\xbar^p =
    (x_1^p,\ldots,x_n^p)$,
  \item \label{iit:hom} or there exist $\phi,\psi \in K[\xbar]$ with
    $\deg_{\xbar}(F) > \max(\deg(\phi), \deg(\psi))$ satisfying the following:
    
    (*) there is an integer $d>1$ \footnote{This condition is actually
      a consequence of $\deg_{\xbar}(F) > \max(\deg(\phi),
      \deg(\psi))$.} and $\ell +1$ polynomials $h_i(u,v) \in K[u,v]$
    homogeneous of degree $d$ such that
$$\left\lbrace
\begin{array}{c}
P(\xbar) = h_0(\phi(\xbar),\psi(\xbar)) = \sum_{k=0}^{d} a_{0k}
\phi(\xbar)^k \psi(\xbar)^{d-k}\\
Q_1(\xbar) = h_1(\phi(\xbar),\psi(\xbar)) = \ldots \hfill \\
\vdots \\
Q_\ell(\xbar) = h_\ell(\phi(\xbar),\psi(\xbar)) = \sum_{k=0}^{d}
a_{\ell k} \phi(\xbar)^k \psi(\xbar)^{d-k}\\
\end{array}
\right.$$
which, setting $H(u,v,\lambdabar) = h_0(u,v) + \sum_{i=1}^\ell
\lambda_i h_i(u,v)$, equivalently rewrites 
$$F(\xbar,\lambdabar) = H(\phi(\xbar),\psi(\xbar),\lambdabar).$$
  \end{enumerate}
\end{enumerate}
\end{theorem}

\begin{remark}
\label{rk:BK}
\begin{enumerate}
  
\item In (\ref{iit:charp}), it follows from $F(\xbar,\lambdabar) \in
  K[\xbar^p,\lambdabar]$ that $P,Q_1,\ldots,Q_\ell$ are in
  $K[x_1^p,\ldots,x_n^p]$; as $K$ is algebraically closed they are
  also $p$-th powers in $K[\xbar]$.
  
\item It follows from the assumption ``$P,Q_1,\ldots,Q_\ell$
  relatively prime'' that the same is true for $\phi$ and $\psi$ in
  (\ref{iit:hom}).
  
\end{enumerate}
\end{remark}

The end of this section is devoted to the study of the decomposition
$F(\xbar,\lambdabar) = H(\phi(\xbar),\psi(\xbar),\lambdabar)$ in
(\ref{iit:hom}) (*) and particularly to the uniqueness of such a
decomposition.

\begin{definition} \label{def:maximal}
  Given two polynomials $\phi,\psi \in K[\xbar]$ relatively prime and such that
  $\deg_{\xbar}(F) > \max (\deg(\phi), \deg(\psi))$,

\begin{enumerate}
\item the polynomial $F$ is said to be
  \defi{$(\phi,\psi)$-homogeneously composed (in degree $d$)} if
  there exists $H(u,v,\lambdabar) \in \overline{K(\lambdabar)}[u,v]$
homogeneous (of degree $d$) in $(u,v)$ such that $F(\xbar,\lambdabar) =
  H(\phi(\xbar),\psi(\xbar),\lambdabar)$.
  The identity $F(\xbar,\lambdabar) =
  H(\phi(\xbar),\psi(\xbar),\lambdabar)$ is then called a
  \defi{$(\phi,\psi)$-homogeneous decomposition of $F$}.
This definition is motivated by  condition (\ref{iit:hom}) (*) of Bertini-Krull theorem.
  
\item A $(\phi,\psi)$-homogeneous decomposition $F(\xbar,\lambdabar) =
  H(\phi(\xbar),\psi(\xbar),\lambdabar)$ is said to be \defi{maximal}
  if $\phi +\lambda \psi$ is irreducible in
  $\overline{K(\lambda)}[\underline x]$ \footnote{where $\lambda$ is a new single variable (to be distinguished from the
    tuple $\lambdabar$).}. 
\end{enumerate}
\end{definition}

\begin{remark}
\begin{enumerate}
\item We also include in this definition the case $\ell = 0$ for which
  only the polynomial $P$ is given. In this situation, the classical
  notion of composed polynomial corresponds to the
  special case of the ``$(\phi,\psi)$-homogeneously composed''
  property for which $\phi$ or $\psi$ is constant.
  
\item For $\ell \geq 1$ we will show that  the \emph{maximality}
  condition is equivalent (except in some special
case) to the maximality of the degree of the
  homogeneous polynomial $H$, whence the terminology.  See theorem
  \ref{thm:unicity} and corollary \ref{cor:degmax}.
  
\item From the Bertini-Krull theorem, ``$\phi +\lambda \psi$
  irreducible in $\overline{K(\underline \lambda)}[\underline x]$'' is
  equivalent to ``$\phi+\lambda\point\psi$ irreducible in $K[\xbar]$
  for at least one $\lambda\point \in K$ with
  $\deg(\phi+\lambda\point \psi)=\max(\deg(\phi),\deg(\psi))$'' and
  also to ``$\phi+\lambda\point\psi$ irreducible in $K[\xbar]$ for
  all but finitely many $\lambda\point \in K$''.
\end{enumerate}
\end{remark}

The
polynomial $F(x,y,\lambda) = x^4-\lambda y^4$ admits the
$(x^2,y^2)$-ho\-mo\-ge\-neous decomposition $F(x,y,\lambda) =
H_1(x^2,y^2,\lambda)$ with $H_1(u,v,\lambda)=u^2-\lambda v^2$. It is
not maximal as $x^2-\lambda y^2 =(x-\sqrt{\lambda}y)
(x+\sqrt{\lambda}y)$. This decomposition however can be refined to a
$(x,y)$-homogeneous decomposition, which is maximal: namely we have
$F(x,y,\lambda) = H_2(x,y,\lambda)$ with $H_2(u,v,\lambda)=u^4-\lambda
v^4$. This refinement is in fact always possible.

\begin{proposition}
\label{prop:maximal}
Assume $F(\xbar,\lambdabar)$ is $(\phi_0,\psi_0)$-homogeneously
composed in degree $d_0$. Then there exists a
maximal $(\phi,\psi)$-homogeneous decomposition of $F$ of degree $d\geq d_0$
and which is of degree $d>d_0$ if  the initial decomposition is not maximal.
\end{proposition}

\begin{proof}
  Let $F(\xbar,\lambdabar) =
  H_0(\phi_0(\xbar),\psi_0(\xbar),\lambdabar)$ be a
  $(\phi_0,\psi_0)$-homogeneous decomposition in degree $d_0$. If
  $\phi_0+\lambda\psi_0$ is irreducible in
  $\overline{K(\lambda)}[\underline x]$ then we are done. Otherwise
  apply the Bertini-Krull theorem to the polynomial
  $\phi_0+\lambda\psi_0$ (note that it is irreducible in
  $K[\lambda][\xbar]$ as $\phi_0$ and $\psi_0$ are relatively prime)
  to conclude that there exist $\phi_1,\psi_1 \in K[\xbar]$ relatively
  prime and with $\max(\deg(\phi_0),\deg(\psi_0)) > \max(\deg(\phi_1),
  \deg(\psi_1))$ such that $\phi_0+\lambda\psi_0$ is
  $(\phi_1,\psi_1)$-homogeneously composed in degree $d_1 \ge 2$.
  Note that this conclusion also covers the extra possibility (\ref{iit:charp})
  of theorem \ref{th:BK} in characteristic $p>0$, which is here  that $\phi_0+\lambda \psi_0$ 
  writes $\phi_1^p + \lambda \psi_1^p$ for
  some $\phi_1,\psi_1\in K[\xbar]$.
  Straightforward calculations on homogeneous polynomials prove that
  $F$ is then $(\phi_1,\psi_1)$-homogeneously composed in degree
  $d_0d_1>d_0$. We can iterate this process, which must stop because
  at each step the degree increases but remains $\leq
  \deg_{\xbar}(F)$. The last step yields a final homogeneous
  decomposition of $F$ which is maximal.
\end{proof}

\subsection{Uniqueness of Bertini-Krull homogeneous
  decompositions}\label{subs:uniqueness}

Theorem \ref{thm:unicity}, which can be viewed as a uniqueness
result for the Bertini-Krull theorem, is the main result of this
section. In this subsection, we assume $\ell \geq 1$.

We need a preliminary adjustment of definition \ref{def:maximal}. 
Given a $(\phi,\psi)$-homogeneous decomposition
$F(\xbar,\lambdabar)=H(\phi,\psi,\lambdabar)$, assume there exists 
$(\alpha, \beta) \not= (0,0)$ in $K^2$ such that
$\alpha \phi + \beta \psi$ is constant in $\underline x$ (that is, is
in $K$).  Then multiplying 
$H(u,v,\lambdabar)$ by any power $(\alpha u + \beta
v)^e$ yields another decomposition $F(\xbar,\lambdabar) = \tilde
H(\phi(\xbar),\psi(\xbar),\lambdabar)$ as above with $\tilde H$
homogeneous (in $u,v$) of degree $\tilde d = d+e$. Conversely if
$H(u,v,\lambdabar)$ has linear factors $\alpha u + \beta v$ (in
$\overline{K(\underline \lambda)}[u,v]$) with $\alpha \phi + \beta
\psi$ constant in $\underline x$, then they are all equal, up to some
constant in $\overline{K(\underline \lambda)}$, to a same linear form
$\alpha_0 u + \beta_0 v \in K[u,v]$ and the homogeneous polynomial
$H\reduced(u,v,\lambdabar)$ obtained from $H(u,v,\lambda)$ by dividing by all possible such
factors $\alpha u + \beta v$ still induces a decomposition
$F(\xbar,\lambdabar) = H\reduced(\phi(\xbar),\psi(\xbar),\lambdabar)$ as
above with $H\reduced$ homogeneous of degree $d\reduced \leq d$. Note we
still have $d\reduced \geq 2$ as $d\reduced \leq 1$ contradicts
$\deg_{\xbar}(F) > \max (\deg(\phi), \deg(\psi))$.

\begin{definition} \label{def:reduced}
  Given two polynomials $\phi,\psi \in K[\xbar]$ relatively prime with
  $\deg_{\xbar}(F) > \max (\deg(\phi), \deg(\psi))$, a
  $(\phi,\psi)$-homogeneous decomposition $F=H(\phi,\psi,\lambdabar)$ is said to be
  \defi{reduced} if the polynomial $H$ has no linear factor $\alpha u
  +\beta v \in K[u,v]$ such that
  $\alpha \phi +\beta \psi$ is constant in $\xbar$.
\end{definition}

From above a reduced $(\phi,\psi)$-ho\-mo\-ge\-neous 
decomposition of $F$ is easily obtained from any $(\phi,\psi)$-homogeneous 
decomposition of $F$.

Also note that if there exists $(\alpha, \beta) \not= (0,0)$ in $K^2$ such that
$\alpha \phi + \beta \psi$ is constant, then up to applying some linear 
transformation $L\in \mathrm{GL}_2(K)$ to $(\phi,\psi)$,  one may assume 
$\phi=1$ and so this can only happen if $F$ is a composed polynomial 
(over $\Klambdabar$). Thus only in this case does
definition \ref{def:reduced} add something to definition \ref{def:maximal}.

\begin{theorem}
\label{thm:unicity}
Assume $\ell \geq 1$. If $F(\xbar,\lambdabar)= P(\xbar) + \lambda_1
Q_1(\xbar) + \cdots + \lambda_\ell Q_\ell(\xbar)$ admits two maximal
homogeneous decompositions:
$$F(\xbar,\lambdabar) =H_1(\phi_1(\xbar),\psi_1(\xbar),\lambdabar)
=H_2(\phi_2(\xbar),\psi_2(\xbar),\lambdabar)$$
then there exists $L\in
\mathrm{GL}_2(K)$ such that $(\phi_1,\psi_1)=L (\phi_2,\psi_2)$.
Furthermore if the two decompositions are reduced then we have $c
\cdot H_2(u,v,\lambdabar)= H_1 (u,v,\lambdabar) \circ L(u,v)$ for some
constant $c\in K$.
\end{theorem}

\begin{example}\label{counter-ex:uniqueness} Theorem \ref{thm:unicity}
  does not extend
  to the case $\ell = 0$. Here is a counter-example. Let $P(x,y)=
  y(x+y)(y^2+xy-2x)$. We have the two maximal 
  homogeneous decompositions:
  
  \noindent - $P=h_1(\phi_1,\psi_1)$ with $h_1(u,v)=v^2-u^2$, $\phi_1=x$,
  $\psi_1=(y-1)(x+y)+y$,  
  
  \noindent - $P=h_2(\phi_2,\psi_2)$ with $h_2(u,v)=uv$, 
  $\phi_2 = y$, $\psi_2=(x+y)(y^2+xy-2x)$.
  
\noindent These two decompositions are distinct even up to elements of
  $\mathrm{GL}_2(K)$.  
\end{example}

\begin{corollary}
\label{cor:degmax} All reduced maximal homogeneous decompositions of $F$
are of the same degree, say $\delta$. Furthermore if $F$ is not a composed 
polynomial over $\Klambdabar$, any homogeneous decomposition of $F$ 
is of degree $\leq \delta$ and equality holds if and only if it is maximal.
\end{corollary}

\begin{proof}[Proof of theorem \ref{thm:unicity}]
  Consider a reduced maximal homogeneous decomposition $F(\underline
  x, \underline \lambda) = H(\phi(\underline x), \psi(\underline
  x),\underline \lambda)$. Write the homogeneous polynomial
  $H(u,v,\underline \lambda)$ (in $u,v$) as a product
  $\prod_{i=1}^d(\alpha_i(\underline \lambda) u+ \beta_i(\underline
  \lambda)v)$ of linear forms in $u,v$ with coefficients in $\overline
  {K(\underline \lambda)}$. Thus we have
  
  $$P(\underline x) + \sum_{i=1}^\ell \lambda_iQ_i(\underline x) =
  \prod_{k=1}^d(\alpha_k(\underline \lambda) \phi(\underline x)+
  \beta_k(\underline \lambda)\psi(\underline x)).$$

\noindent
The result will be easily deduced from these two claims and the unique
factorization property in the domain $\overline {K(\underline
  \lambda)}[\underline x]$.
\begin{itemize}
\item[(a)] There are at least two factors $\alpha_k(\underline
  \lambda) \phi(\underline x)+ \beta_k(\underline
  \lambda)\psi(\underline x)$ that are non constant in $\xbar$ and non
  proportional (by some constant in $\Klambdabar$).
  
\item[(b)] All factors $\alpha_k(\underline \lambda) \phi(\underline
  x)+ \beta_k(\underline \lambda)\psi(\underline x)$ ($k=1,\ldots,d$)
  are irreducible in $\overline {K(\underline \lambda)}[\underline x]$
  and are not in $K[\underline x]$ (even up to constants in
  $\Klambdabar$).
\end{itemize}

{\it Proof of claim (a).} First note that due to definition
\ref{def:reduced}, no factor $\alpha_k \phi+ \beta_k\psi$ is in
$\overline {K(\underline \lambda)}$. Assume (a) does not hold. Then
$F(\underline x, \underline \lambda)$ is of the form $\alpha M^d$ with
$\alpha \in \overline {K(\underline \lambda)}$ and $M=\alpha_1 \phi +
\beta_1\psi$.
Taking the derivative with respect to $\lambda_i$ shows that $M^{d-1}$
divides $Q_i$ in $\overline {K(\underline \lambda)}[\underline x]$,
$i=1,\ldots,\ell$. But as $M^d$ divides $F(\underline x, \underline
\lambda)$, we obtain that $M^{d-1}$ divides $P$ as well.
A contradiction as $\deg(M)>0$ and $P,Q_1,\ldots,Q_\ell$ are assumed
to be relatively prime.

{\it Proof of claim (b).} Assume that for some $k\in \{1,\ldots,d\}$,
$\alpha_k \phi + \beta_k\psi$ is reducible in $\overline {K(\underline
  \lambda)}[\underline x]$. One may assume that $\deg(\psi) >0$ and
$\beta_k\not=0$.  If $\alpha_k \not= 0$, set $\mu(\underline \lambda)
= \beta_k(\underline \lambda)/ \alpha_k(\underline \lambda)$. The
polynomial $ \phi(\underline x) + \mu(\underline \lambda)
\psi(\underline x)$ is reducible in $\Klambdabar[\xbar]$ and
consequently so are the polynomials $ \phi(\underline x) +
\mu(\underline \lambda\point) \psi(\underline x)$ for all
specializations $\underline \lambda \rightarrow \underline
\lambda\point$ in $K^\ell$ except possibly in a proper
Zariski closed subset.  
It follows then from the Bertini-Krull theorem and the
irreducibility of $\phi+\lambda \psi$ in $\overline
{K(\lambda)}[\underline x]$ that $ \mu(\underline \lambda)$ has only
finitely many specializations in $K$ and so necessarily
$\mu(\underline \lambda)=\mu \in K$. Then set $a(\underline
x)=\phi(\underline x) + \mu \psi(\underline x)$. In the case that
$\alpha_k = 0$, set $a(\underline x) = \psi(\underline x)$. In all
cases, $a(\underline x)\in K[\underline x]\setminus K$ and
$F(\underline x,\underline \lambda) = a(\underline x) G_{\underline
  \lambda}(\underline x)$ for some $G_{\underline \lambda}(\underline
x)\in \overline {K(\underline \lambda)}[\underline x]$. We now show
that this leads to a contradiction. Namely for each $i=1,\ldots,\ell$
$$\frac{\partial G_{\underline \lambda}}{\partial \lambda_i} =
\frac{1}{a} \frac{\partial F}{\partial \lambda_i}=\frac{Q_i}{a}$$
lies
both in $K(\underline x)$ and in $\overline {K(\underline
  \lambda)}[\underline x]$, and so is in $K[\underline x]$. Thus $a$
divides $Q_i$ in $K[\underline x]$, $i=1,\ldots,\ell$. But as $a$
divides $P+\sum_{i=1}^\ell \lambda_i Q_i$, $a$ divides $P$ as well 
(both in $\Klambdabar[\xbar]$): a
contradiction as $\deg(a)>0$ and $P,Q_1,\ldots,Q_\ell$ are relatively
prime.

It follows from claims (a) and (b) that if the two maximal homogeneous
decompositions given in the statement of theorem \ref{thm:unicity}
are reduced, then we have $(\phi_1,\psi_1) = L_\lambdabar
(\phi_2,\psi_2)$ for some $L_\lambdabar \in
\mathrm{GL}_2(\overline{K(\underline \lambda)})$. Now for all
$\underline \lambda\point \in K^\ell$ but in a proper Zariski closed
subset we also have $(\phi_1,\psi_1) = L_{\lambdabar\point}
(\phi_2,\psi_2)$ with $L_{\lambdabar\point} \in \mathrm{
  GL}_2(K)$.

It also follows from claims (a) and (b) that the set of linear factors
$\alpha_k(\underline \lambda) u+ \beta_k(\underline \lambda)v $ of the
polynomial $H(u,v,\underline \lambda)$ is uniquely determined (up to
non zero constants) by the set of irreducible factors
$\alpha_k(\underline \lambda) \phi(\underline x)+ \beta_k(\underline
\lambda)\psi(\underline x)$ of $F(\underline x,\underline \lambda)$.
This yields the additional conclusion $c \cdot H_2(u,v,\lambdabar)=
H_1 (u,v,\lambdabar) \circ L_{\lambdabar\point}(u,v)$ of theorem \ref{thm:unicity}.

Finally if the two given maximal homogeneous decompositions of $F$ are
not reduced, consider the two associated reduced decompositions $F
=H_1\reduced(\phi_1,\psi_1,\lambdabar)
=H_2\reduced(\phi_2,\psi_2,\lambdabar)$ (constructed prior to definition
\ref{def:reduced}). The proof above still yields $(\phi_1,\psi_1)
= L_{\lambdabar\point} (\phi_2,\psi_2)$ for some $L_{\lambdabar\point} \in \mathrm{GL}_2(K)$.
\end{proof}

\subsection{Further comments}
\label{subsection:comments}
Retain the notation from the above proof.

\subsubsection{} As a consequence of the factors $\alpha_k(\underline
\lambda) \phi(\underline x)+ \beta_k(\underline
\lambda)\psi(\underline x)$ not being in $K[\underline x]$ even up to
constants in $\Klambdabar$ we have $\alpha_k(\underline \lambda)
\beta_k(\underline \lambda) \not=0$ and $\deg_{\underline
  x}(\alpha_k\phi+ \beta_k\psi) = \max(\deg(\phi),\deg(\psi))$,
$k=1,\ldots,d$.

\subsubsection{} From the Bertini-Noether theorem  \cite[proposition 8.8]{FrJa},
for all $\underline \lambda\point \in K^\ell$ but in a proper Zariski
closed subset $\mathcal Z$, the polynomials $\alpha_k(\underline
\lambda\point) \phi(\underline x)+ \beta_k(\underline
\lambda\point)\psi(\underline x)$, obtained by specializing
$\underline \lambda$ to $\underline \lambda\point$ in the irreducible
factors $\alpha_k(\underline \lambda) \phi(\underline x)+
\beta_k(\underline \lambda)\psi(\underline x)$ of $F(\underline
x,\underline \lambda)$, are the irreducible factors of $F(\underline
x,\underline \lambda\point)$ in $K[\underline x]$.

\subsubsection{} The vector space $\overline{K(\underline \lambda)} \phi +
\overline{K(\underline \lambda)} \psi$, which is uniquely determined
by $F(\underline x,\underline \lambda)$, is the
$\overline{K(\underline \lambda)}$-vector space generated by all
irreducible divisors of $F(\underline x,\underline \lambda)$ in
$\overline{K(\underline \lambda)}[\underline x]$. As to the $K$-vector
space $K\phi + K\psi$, it is the vector space generated by all
irreducible divisors in $K[\underline x]$ of the polynomials
$F(\underline x,\underline \lambda\point)$ with $\underline
\lambda\point\notin \mathcal Z$ (where $\mathcal Z$ is defined just
above).

\subsubsection{} Consider the problem, given a polynomial $P$ as above, of
finding all the sets $\{Q_1,\ldots,Q_\ell\}$ of polynomials as above
(with $\ell\geq 1$), such that $P+\lambda_1 Q_1+\cdots+\lambda_\ell
Q_\ell$ is reducible in $\Klambdabar[\xbar]$. This problem will be
studied in the next section in the special situation
$Q_1,\ldots,Q_\ell$ are monomials. We note here that the general
problem can be reduced to the special case $\ell=1$. 

Indeed, if
$\{Q_1,\ldots,Q_\ell\}$ is a solution to this problem, then,
for some integer $d\geq 2$,
the polynomials $P,Q_1,\ldots, Q_\ell$ all are in the $d$-th symmetric
power $(K\phi + K\psi)^d$
of some vector space $K\phi + K\psi \subset K[\xbar]$ \footnote{This is another 
way of saying that each of these
  polynomials can be written $h(\phi,\psi)$ with $h\in K[u,v]$
  homogeneous of degree $d$.} which from theorem \ref{thm:unicity}
   is uniquely determined by
$P,Q_1,\ldots,Q_\ell$. Now there exists $Q\in (K\phi + K\psi)^d$ that
is relatively prime to $P$. Clearly $P+\lambda Q$ is reducible in
$\Klambdar[\xbar]$, that is, the singleton $\{Q\}$ is a solution to
the problem with $\ell=1$.  The vector space $K\phi + K\psi$ is also
uniquely determined by $P$ and $Q$.  Thus finding all solutions $Q$ to the
problem with $\ell=1$ provides all possible solutions
$\{Q_1,\ldots,Q_\ell\}$ to the general problem: these sets are all
possible finite subsets of the sets $(K\phi + K\psi)^d$ attached to
the solutions $Q$. 

For self-containedness of next section, we
will not use this remark there. We just state this other related consequence 
of theorem \ref{thm:unicity}.

\begin{corollary}
Suppose given two maximal homogeneous decompositions
$P(\xbar)+\lambda_1 Q_{1}(\xbar)+\cdots+\lambda_{\ell} Q_{\ell}(\xbar)=H(\phi(\xbar),\psi(\xbar),\lambdabar)$
and $P(\xbar)+\lambda^\prime_1 Q^\prime_{1}(\xbar)+\cdots+\lambda^\prime_{\ell^\prime} Q^\prime_{\ell^\prime}(\xbar)=H^\prime(\phi^\prime(\xbar),\psi^\prime(\xbar),\lambdabar^\prime)$ (with $\ell, \ell^\prime \geq 1$). Assume further  that $Q_1=Q^\prime_1$ and that
$P$ and $Q_1$ are relatively prime. Then we have $(\phi',\psi') = L(\phi,\psi)$
  for some $L \in \mathrm{GL}_2(K)$.
\end{corollary}

\section{Reducibility monomial sites}
\label{sec:mon.sites}

We keep the notation of section \ref{sec:uniqueness} but assume in
addition that $\ell\geq 1$ and that $Q_1,\ldots,Q_\ell$ are monomials
such that $\deg(Q_i)\leq \deg(P)$, $i=1,\ldots,\ell$. We set
$Q_i=x_1^{e_{i1}} \cdots x_n^{e_{in}}$, $i=1,\ldots,\ell$.

\begin{definition}
  The set $\{Q_1,\ldots,Q_\ell\}$ is said to be a \defi{reducibility
    monomial site} of $P$ is
  $F(\xbar,\lambdabar)=P+\lambda_1Q_1+\cdots+\lambda_\ell Q_\ell$ is
  reducible in $\Klambdabar[\xbar]$. If $\ell=1$ we just say $Q_1$ is a \defi{reducibility
  monomial}.
\end{definition}

It is readily checked that any subset of a reducibility monomial site
is a reducibility monomial site.

\begin{definition}
  A polynomial $P\in K[\xbar]$ is said to
  be \defi{homogeneous in two monomials} if $P$ is
  $(m_1,m_2)$-homogeneously composed for some monomials
  $m_1$ and $m_2$ (which according to definition \ref{def:maximal}
  should be relatively prime and such that $\deg(P) > \max(\deg(m_1),
  \deg(m_2))$).
\end{definition}

This property can be easily detected thanks to the Newton  
representation of $P$ (as already used in the introduction).
Indeed, set $m_1=x_1^{a_1}\cdots x_n^{a_n}$ and $m_2=x_1^{b_1}\cdots
x_n^{b_n}$. If $P$ is homogeneous in $m_1$ and $m_2$, then $P$ is a sum of 
monomials of the form:
$$m_1^k m_2^{d-k} = x_1^{db_1+k(a_1-b_1)} \cdots x_n^{db_n+k(a_n-b_n)}
\hskip 5mm (k\in \{0,\ldots,d\})$$

\noindent The corresponding points $M_k=(db_1+k(a_1-b_1),
\ldots,db_n+k(a_n-b_n))$ ($k=0,\ldots,d$)  lie on a straight line in $\Qq^n$.\footnote{Note 
however that the monomials being lined up in the Newton representation is not sufficient 
for $P$ to be  homogeneous in two monomials: for example $P=xy+x^2y^4 + x^3y^6$
has that property but is not homogeneous in two monomials. It is of course easy to give
a full test for some polynomial $P$ to be homogeneous in terms of its Newton 
representation but writing out the exact condition is not very enlightening. See
also remark \ref{remark:newton}. }

We will show below (theorem \ref{th:mon} (addendum 1)) that a $(m_1,m_2)$-homogeneous
  decomposition of $P$ is maximal, that is $m_1+\lambda m_2$ is
  irreducible in $\Klambdar[\xbar]$ if and only if $m_1$ and $m_2$
  are not $d$-th powers in $K[\xbar]$
  for some integer $d>1$, or,
  equivalently, if $a_1, \ldots, a_n, b_1,\ldots,b_n$ are
  relatively prime.


\subsection{Main theorem}
Our main result determines the reducibility monomial sites of a polynomial. We first state it
in the general situation of a polynomial that is neither a mo\-no\-mial nor a pure power. 
The two remaining special cases are dealt with in two addenda. The proof 
is given in section \ref{proofMT}.

\begin{theorem}[general case]
\label{th:mon}
Assume $P(\xbar)$ is not a monomial and is not a pure power in
$K[\xbar]$.
\begin{enumerate}
\item \label{it:hom} If $P$ is homogeneous in two monomials, then given 
  a maximal $(m_1,m_2)$-homogeneous decomposition $P=h(m_1,m_2)$ of
  degree $\delta$ with $m_1$ and $m_2$ monomials\footnote{Such a
    decomposition exists (proposition \ref{prop:maximal}) and is
    unique up to trivial transformations (lemma
    \ref{lem:monunicity}).}, the reducibility monomial sites of $P$ are
  all sets of monomials $m_1^km_2^{\delta-k}$, $0 \le k \le
  \delta$, of degree $\leq \deg(P)$.
\smallskip
  
\item \label{it:mon} If $P$ is not homogeneous in two monomials 
 then the only possible reducibility monomial sites are singletons 
  ($\ell=1$) of the form $\{m^d\}$ with $m$ a monomial relatively prime
  to $P$ and $d\ge 2$. Furthermore the following should hold: $P = h(m,\psi)$ 
  with $h\in K[u,v]$ homogeneous of degree $d$, $\psi \in K[\xbar]$ non monomial
  and $\deg(P) > \max(\deg(m),\deg(\psi))$\footnote{By 
proposition \ref{prop:maximal} we may also impose that $\psi + 
\lambda m$ is irreducible in $\Klambdar[\xbar]$.}.
\end{enumerate}
\end{theorem}


\begin{remark}
\begin{enumerate}
\item In the homogeneous case (1), the reducibility monomials
  $m_1^km_2^{\delta-k}$ also are on the line formed by the monomials
  of $P$ in its Newton representation.
  
\item In case (\ref{it:mon}) we do not know whether there may be
  several reducibility monomials of the form $m^d$. This is related to
  the possibility that $P$ can be written $P = h(m,\psi)$ as in the
  statement in several different ways, and so to the uniqueness of
  homogeneous decompositions of $P$. In section \ref{subs:uniqueness}
  where this problem is studied for the polynomial
  $P+\lambda_1Q_1+\cdots + \lambda_\ell Q_\ell$ with $\ell \geq 1$, we
  give a counter-example to uniqueness for $\ell = 0$ (example
  \ref{counter-ex:uniqueness}).  However the two monomials $m^d$
  associated to the two homogeneous decompositions of $P$ shown there
  are $x^2$ and $y^2$; the second one is not relatively prime to $P$
  and so is not a reducibility monomial according to our definitions.
  
\item In case (2) where $P=h(m,\psi)$, by setting $g(t)= h(1,t)$ we obtain $P/m^d = g(\psi/m)$
is a \textit{composite rational function}
  as considered in \cite{Bo} (of special form though as $g$ is here
  a polynomial).
\end{enumerate}
\end{remark}

\subsection{The monomial case}
\label{monomial case}

Here we consider the case $P$ is a monomial $\gamma \hskip
2pt x_1^{e_{1}} \cdots x_n^{e_{n}}$ (with $\gamma \in K,
\gamma\not=0$). The argument below can be viewed as an easy
special case of the general method.

From \S \ref{sec:uniqueness}, if $F(\xbar,\lambdabar)$ is reducible in
$\Klambdabar[\xbar]$, then equivalently either $F(\xbar,\lambdabar) \in
K[\xbar^p,\lambdabar]$ (with $\charac(K)=p>0$) or $F(\xbar,\lambdabar)$ is
$(\phi,\psi)$-homogeneously composed in degree $d$ for some
$(\phi,\psi) \in K[\xbar]$.  In the latter case, factor
the homogeneous polynomials involved in the decomposition as products
of linear forms to obtain

$$\left\lbrace
\begin{array}{c}
P(\xbar) = \prod_{k=1}^{\mu_0} (\alpha_{0k} \phi(\xbar) +
\beta_{0k}\psi(\xbar))^{r_{0k}} \hfill\\ 
Q_i(\xbar) = \prod_{k=1}^{\mu_i} (\alpha_{ik} \phi(\xbar) +
\beta_{ik}\psi(\xbar))^{r_{ik}}\hskip 5mm (i=1,\ldots,\ell)\\
\end{array}
\right.$$

\noindent
where the $(\alpha_{ik},\beta_{ik})$ are non-zero, pairwise non proportional and the integers $r_{ik}$ are $>0$ and satisfy
$\sum_{k=1}^{\mu_i} r_{ik} = d$ ($i=0,\ldots,\ell$).

All the factors appearing in the right-hand side terms are necessarily
monomials and at least two of them are non proportional (as
$P,Q_1,\ldots,Q_\ell$ are relatively prime). Therefore up to changing
$(\phi,\psi)$ to $L(\phi,\psi)$ for some $L\in \mathrm{GL}_2(K)$ one
may assume that $\phi$ and $\psi$ themselves are two monomials $m_1$
and $m_2$. Taking into account that $P,Q_1,\ldots,Q_\ell$ are
monomials and that they are relatively prime, we obtain the following
characterization (the converse is clear).

\begin{addendumA}[addendum 1]
If $P$ is a monomial the following are equivalent:
\begin{enumerate}
  
\item The polynomial $P+\lambda_1Q_1+\cdots+\lambda_\ell Q_\ell$ is
  reducible in $\Klambdabar[\xbar]$ (that is, $\{Q_1,\ldots,Q_\ell\}$ is a reducibility
    monomial site of $P$),
  
\item
  \begin{enumerate}
  \item either $\charac K = p >0$ and $P,Q_1,\ldots,Q_\ell \in
    K[\xbar^p]$,
  \item or $P,Q_1,\ldots,Q_\ell$ are of the form $m_1^k m_2^{d-k}$
    ($0\leq k\leq d$) for some relatively prime monomials $m_1$ and
    $m_2$ and some integer $d>1$, and they include $m_1^d$ and
    $m_2^d$.
    
    Furthermore, for all $(\phi,\psi)$-homogeneous decompositions of
    $P+\lambda_1Q_1+\cdots+\lambda_\ell Q_\ell$, $(\phi,\psi)$ is a
    couple of monomials, up to some element $L\in \mathrm{GL}_2(K)$.
  \end{enumerate}
\end{enumerate}
\end{addendumA}

\begin{remark}
  In general there may be several couples $(m_1,m_2)$ such that $P$ is
  of the form $m_1^k m_2^{d-k}$,
    and so
  several corresponding reducibility sites for $P$. For example
  $P=x^3y^2$ is homogeneously composed for both couples of
  monomials $(x^3,y^2)$ and $(x^3y^2,1)$ and both decompositions are maximal.
  In the non monomial case, this will not happen: up to trivial
  transformations the couple $(m_1,m_2)$ is uniquely determined by $P$
  (see lemma \ref{lem:monunicity}).
\end{remark}

\subsection{Pure power case} \label{subsection:purepower}In the case 
$P$ is a pure power in $K[\xbar]$, the three
following possibilities can occur:

\begin{enumerate}
\item $P$ is homogeneous in two monomials. In this case let $P = h(m_1,m_2)$ be a maximal homogeneous decomposition
of degree  $\delta$ in two monomials $m_1$ and $m_2$ and set $\mathcal{M}_1 = \{ m_1^km_2^{\delta-k} \mid 0 \le k \le \delta\}$. All subsets of $\mathcal{M}_1$
are reducible monomial sites.
  
\item $P$ admits a maximal $(m,\psi)$-homogeneous decomposition in degree $d$,
 with $m$ a monomial and $\psi\in K[\xbar]$ non monomial. In this case, if $\deg(m^d) \leq \deg(P)$,  then $m^d$ is a reducibility monomial.

\item $\charac(K)=p>0$ and $P\in K[\xbar^p]$. In this case set $\mathcal{M}_3 = \{
  m^p \mid m \text{ is a monomial and} \deg(m^p) \leq \deg (P)\}$. All subsets of $\mathcal{M}_3$
  are reducible monomial sites.
\end{enumerate}

\begin{addendumB}[addendum 2]
Assume $P$ is a pure power but is not a monomial. Then the reducibility monomial sites of $P$ are those described in possibilities (1), (2) and (3).
\end{addendumB}

The following observations make the pure power case rather special:
\vskip 1mm

(a) possibility (2) is always satisfied: indeed by assumption we have $P=S^e$ for some $S\in K[\xbar]$ and some integer $e>1$, which is a $(m,S)$-homogeneous decomposition of degree $e$ for any monomial $m$
relatively prime to $S$; the corresponding monomials $m^e$ with $\deg(m^e)\leq \deg(P)$ are reducibility monomials. However there may be other kinds of decompositions $P=h(m,\psi)$. For example, take $P(x,y) = (2y^3-x^4)^2x^4$. Squares monomials of degree $\leq 12$ are reducibility monomials. Now for $m = y^3$, $\psi = y^3-x^4$ and $h(u,v)=(u+v)^2(u-v)$, we also have $P= h(m,\psi)$ and so $m^3=y^9$ is another reducibility monomial of $P$.
\vskip 1mm

(b) possibilities (1), (2) and (3) can occur simultaneously. Take for example $P(x,y) = (x^2-y^3)^3$.  Then $P$ is homogeneous in the two monomials $x^2$ and $y^3$; the corresponding set $ \mathcal{M}_1$ is $\mathcal{M}_1=\{x^6, x^4y^3, x^2y^6, y^9\}$.  As $P$ is a third power, each of the monomials $1$, $x^3$, $y^3$, $x^6$,  $x^3y^3$, $y^6$, $x^9$, $x^6y^3$, $x^3y^6$, $y^9$ is a reducibility monomial. Finally if $\charac (K) = 3$, then every subset of $\mathcal{M}_3=\{1,x^3, y^3, x^6,  x^3y^3, y^6, x^9, x^6y^3, x^3y^6, y^9\}$ is a reducibility monomial site.

\subsection{Lemmas}

The two following lemmas will be used in the proof of theorem \ref{th:mon}.

\begin{lemma}
\label{lem:monirr}
Given two monomials $m_1,m_2 \in K[\xbar]$ such that we have
$\max(\deg(m_1), \deg(m_2))>0$, the following are equivalent:

\begin{enumerate}
\item[(i)] there exists $\lambda\point \in K$, $\lambda\point \neq
  0$, such that $m_1+\lambda\point m_2$ is irreducible in $K[\xbar]$,
  
\item[(ii)] for all $\lambda\point \in K$, $\lambda\point\not=0$,
  $m_1+\lambda\point m_2$ is irreducible in $K[\xbar]$,
  
\item[(iii)] $m_1+\lambda m_2$ is irreducible in $\Klambdar[\xbar]$.
\end{enumerate}
\end{lemma}

\begin{proof} The equivalence (iii)$\Leftrightarrow$(i) is a special case
  of the Bertini-Krull theorem and (ii)$\Rightarrow$(i) is trivial. We
  are left with proving (i)$\Rightarrow$(ii). Assume there exist
  $\lambda\point_1, \lambda\point_2 \in K$, both non zero and
  such that $m_1+\lambda\point_1 m_2$ is reducible and
  $m_1+\lambda\point_2 m_2$ is irreducible in $K[\xbar]$.
  
  Set $m_1 = x_1^{a_1} \cdots x_n^{a_n}$ and $m_2 = x_1^{b_1} \cdots
  x_n^{b_n}$. One may assume that $\deg(m_2)>0$ and so for example
  $b_1>0$.  If $a_1>0$ then $x_1$ divides
  $m_1+\lambda\point_2 m_2$ and so $m_1=m_2=x_1$ (up to some non-zero 
  multiplicative constants) in which case the result is obvious. 
  Thus one may assume $a_1=0$. If $m_1(\xbar)+\lambda\point_1 m_2(\xbar) =
  R(\xbar) \cdot S(\xbar)$ is a non trivial factorization of
  $m_1+\lambda\point_1 m_2$ ($\deg(R), \deg(S) > 0$), we have
  $$m_1+\lambda\point_2 m_2 = R\left(({\lambda\point_1}^{-1}\lambda\point_2
    )^\frac{1}{b_1}x_1,x_2,\ldots,x_n\right)\cdot
  S\left(({\lambda\point_1}^{-1}\lambda\point_2
    )^\frac{1}{b_1}x_1,x_2,\ldots,x_n\right)$$
  which contradicts the irreducibility of $m_1+\lambda\point_2
  m_2$.
\end{proof}

\begin{lemma}
\label{lem:monunicity}
Assume $P(\xbar)$ is not a monomial and is given with a maximal
$(m_1,m_2)$-homogeneous decomposition $P=h(m_1,m_2)$ of degree $d$
with $m_1$ and $m_2$ monomials.

\begin{enumerate}
\item \label{it:monuni1} If $P = h^\prime(m^\prime_1,m^\prime_2)$ is another maximal
  homogeneous decomposition of degree $d^\prime$ of $P$ in monomials
  $m^\prime_1$ and $m^\prime_2$, then either ($m_1=a m^\prime_1$ and
  $m_2=b m^\prime_2$) or ($m_1=a m^\prime_2$ and $m_2=b m^\prime_1$),
  for some non-zero constants $a,b\in K$, and $d=d^\prime$.
  
\item \label{it:monuni2} There is no maximal homogeneous $(m,\psi)$-decomposition of $P =
  h^\prime(m,\psi)$ with $\psi \in K[\xbar]$ non monomial and 
$m$ a monomial relatively prime to $P$ and not a monomial of $\psi$
 unless $P=\psi^{d''}$ with $\psi$ homogeneous in
  $m_1$ and $m_2$ and $d'' \geq 2$.
\end{enumerate}
\end{lemma}

\begin{proof}
  We can write
  $$P= h(m_1,m_2)=\prod_{k=1}^{\mu} (\alpha_k m_1 + \beta_k m_2)^{r_k}
  \leqno{(**)}$$
  where the $(\alpha_k,\beta_k)$ are non-zero and
  pairwise non-proportional and the integers $r_k$ are $>0$ and satisfy
  $\sum_{k=1}^\mu r_k = d$.
  
  \smallskip
  
  (1) As $P$ is not a monomial there exists $k \in \{1,\ldots,\mu\}$
  such that $\alpha_k\beta_k \neq 0$.  Then by lemma \ref{lem:monirr}
  $\alpha_k m_1 + \beta_k m_2$ is irreducible in $K[\xbar]$.
  
  Assume $P$ has another maximal homogeneous decomposition in
  monomials $m^\prime_1$ and $m^\prime_2$
  
  $$P = h^\prime(m^\prime_1,m^\prime_2)=\prod_{k=1}^{\mu^\prime}
  (\alpha_k^\prime m_1^\prime + \beta'_k m_2^\prime)^{r'_k}$$

\noindent
where the $(\alpha_k^\prime,\beta_k^\prime)$ are non-zero, pairwise non proportional
 and the integers $r_k^\prime$ are $>0$ and satisfy
$\sum_{k=1}^{\mu^\prime} r'_k >1$. From the unique factorization
property in the domain $K[\xbar]$, there exists $h \in
\{1,\ldots,\mu^\prime\}$ with $\alpha_h^\prime \beta_h^\prime \neq 0$
such that, up to a non-zero multiplicative constant, we have
$\alpha_k m_1 + \beta_k m_2 = \alpha_h^\prime m_1^\prime + \beta'_h
m_2^\prime$.
As $m_1, m_2, m^\prime_1, m^\prime_2$ are monomials we obtain the
desired conclusion.  

\smallskip

\begin{remark}
\label{remark:newton}
 In fact the monomials $m_1$ and $m_2$ of some maximal homogeneous decomposition
  of $P$ can be easily recovered from the Newton representation of $P$.
  Indeed, using the notation from the beginning of section \ref{sec:mon.sites}, for any 
  two distinct points $M_h$ and $M_k$,
  we have $\overrightarrow{M_kM_h} = (k-h) \overrightarrow \Delta$ where
  $\overrightarrow \Delta = (a_1-b_1,\ldots, a_n-b_n)$. As $\min(a_j,b_j)=0$, $j=1,\ldots,\ell$, the non-zero exponents
  of $m_1$ (resp. of $m_2$) correspond to the positive components
  (resp. to the negative components) of $\overrightarrow \Delta$. As
  $a_1, \ldots, a_n,b_1,\ldots,b_n$ are relatively prime, these
  exponents correspond to the components of $\overrightarrow{M_kM_h}$
  divided by their g.c.d. 
  
\end{remark}

(2) Suppose $P$ has a maximal $(m,\psi)$-homogeneous decomposition
(with $m$ and $\psi$ as in the statement)

$$P = h^\prime(m,\psi) = \prod_{k=1}^{\mu^\prime} (\alpha_k^\prime
\psi + \beta_k^\prime m)^{r_k^\prime}$$

\noindent
where the $(\alpha_k^\prime,\beta_k^\prime)$ are non-zero and pairwise non-proportional 
 and the integers $r_k^\prime$ are $>0$ and satisfy
$\sum_{k=1}^{\mu^\prime} r_k^\prime = d^{\prime}>1$.

Consider first the case there exists $h\in \{1,\ldots,\mu^\prime\}$
with $\alpha_h^\prime \beta_h^\prime \not= 0$. Comparing with $(**)$
above we obtain that the polynomial $\alpha_h^\prime \psi +
\beta_h^\prime m$ is a product of say $\nu$ irreducible factors
$\alpha_k m_1 + \beta_k m_2$ with $\alpha_k \beta_k \not = 0$ 
(irreducible by lemma \ref{lem:monirr})
and possibly some monomial $\rho$. As
$\alpha_h^\prime \psi + \beta_h^\prime m$ has at least $3$ monomials,
the integer $\nu$ is $\geq 2$. Thus $\alpha_h^\prime \psi +
\beta_h^\prime m$ can be written $\rho \hskip 2pt \kappa(m_1,m_2)$
with $\kappa\in K[u,v]$ homogeneous of degree $\nu \geq 2$. As $m$ is not a monomial of $\psi$, conclude
that, up to non zero constants in $K$, $m$ is one of the monomials of
$\kappa(m_1,m_2)$ multiplied by $\rho$ and that $\psi$ is the sum of
the other monomials of $\kappa(m_1,m_2)$, also multiplied by $\rho$.
Now as $\psi$ and $m$ are relatively prime, $\rho$ is a non-zero
constant in $K$. But then $m+\lambda \psi$ is $(m_1,m_2)$-homogeneously
composed in degree $\nu$, which contradicts the maximality of the
$(m,\psi)$-decomposition.

Assume next that $\alpha_h^\prime \beta_h^\prime = 0$ for all
$h=1,\ldots, \mu^\prime$. If no coefficient $\alpha_h^\prime$ is zero,
then $P=\psi^{d^\prime}$ (up to some non-zero multiplicative constant).  
If some coefficient $\alpha_h^\prime$ is
zero, then $m$ divides $P$ and as $P$ and $m$ are assumed to be
relatively prime, $m$ is a non-zero constant in $K$. Conclude in both
cases that $P=\psi^{d''}$ (up to some non-zero multiplicative constant)
where $d''$ is the number of coefficients
$\alpha_h^\prime$ that are non-zero (counted with the multiplicities $r^\prime_k$); 
we have $d'' \leq d^\prime$ and
$d''\geq 2$ for otherwise we would have $\deg(P)\leq \max(\deg(\psi),
\deg(m))$.  Observe next that the exponents $r_{k}$ are all divisible
by $d''$: if $\alpha_k \beta_k\not= 0$, this is because $\alpha_km_1+
\beta_km_2$ is irreducible in $K[\xbar]$ and for the possible two
factors that are powers of $m_1$ and $m_2$, because $m_1$ and $m_2$
are relatively prime.  Conclude that $\psi$ is as announced
homogeneous in $m_1$ and $m_2$.
\end{proof}

\subsection{Proof of theorem \ref{th:mon}} 
\label{proofMT} Addendum 1 has already been proved (in section \ref{monomial case}) so
we may assume $P$ is not a monomial.

\subsubsection{Preliminary discussion:} Let $\{Q_1,\ldots, Q_\ell\}$ ($\ell
\ge 1$) be a reducibility monomial site of $P$.  

From remark \ref{rk:BK} the case (\ref{iit:charp}) in
the Bertini-Krull theorem can only occur if $P$ is a pure power, and in this
case the conclusion corresponds to possibility (3) of theorem \ref{th:mon} (addendum 2).

Suppose now it is part (\ref{iit:hom}) of the Bertini-Krull theorem that holds. That is, 
the polynomial $F(\xbar,\lambdabar) =
P+\lambda_1Q_1+\cdots+\lambda_\ell Q_\ell$ has a
$(\phi,\psi)$-homogeneous decomposition in degree $d$ for some
$\phi, \psi\in K[\xbar]$, which in addition we may and will assume
to be maximal (proposition \ref{prop:maximal}).

Thus we have $P(\xbar)= h_0(\phi(\xbar),\psi(\xbar))$ and $Q_i(\xbar)=
h_i(\phi(\xbar),\psi(\xbar))$ ($i=1,\ldots,\ell$) for some homogeneous
polynomials $h_0,\ldots, h_\ell \in K[u,v]$ of degree $d$. Note that
as $\deg(P)\geq \deg(Q_i)$, $i=1,\ldots,\ell$, we have
$\deg_{\xbar}(F) = \deg(P) > \max(\deg(\phi),\deg(\psi))$ and so
$P=h_0(\phi,\psi)$ is still a $(\phi,\psi)$-homogeneous decomposition
of $P$.  Write then $h_i(u,v) = \prod_{k=1}^{\mu_{i}} (\alpha_{ik} u +
\beta_{ik} v)^{r_{ik}}$ with, for each $i=0,\ldots,\ell$, the
$(\alpha_{ik},\beta_{ik})$ non-zero and pairwise non proportional and the
integers $r_{ik}>0$ and satisfying $\sum_{k=1}^{\mu_i} r_{ik}=d$.
Unless $\ell = 1$ and $Q_1$ is constant, one may assume $Q_1$ is a non
constant monomial and then all factors $\alpha_{1k} \phi(\xbar) +
\beta_{1k} \psi(\xbar)$ ($k=1,\ldots,\mu_1$) are monomials and at
least one, say $m$, is non constant.  If $\ell = 1$ and $Q_1$ is
constant, then $\phi$ or $\psi$, say $\phi$ is constant. In all cases,
up to changing $(\phi,\psi)$ to $L(\phi,\psi)$ for some $L \in
\mathrm{GL}_2(K)$, one may assume that $\phi$ is a monomial $m$
and that $m$ is not a monomial of $\psi$.
Observe then that if $\psi$ has at least two monomials 
then $Q_i=h_i(m,\psi)$ can be a monomial only if $h_i(u,v)=u^d$
and so $\ell = 1$ and $Q_1=m^d$.

\medskip We now distinguish two cases.

\subsubsection{1st case:} $P$ is homogeneous in two monomials.

Let $P=h(m_1,m_2)$ be a maximal $(m_1,m_2)$-homogeneous decomposition
in degree $\delta$ with $m_1$ and $m_2$ monomials. From above
$P=h_0(m,\psi)$ is another maximal homogeneous decomposition.

If $\psi$ itself is a monomial then from lemma
\ref{lem:monunicity} (\ref{it:monuni1}), we have $d=\delta$ and $(m,\psi)=(am_1,bm_2)$ or
$(m,\psi)=(bm_2,am_1)$ for some non-zero constants $a,b\in K$.
Conclude each $Q_i$ is homogeneous in $m_1$ and $m_2$ in degree
$\delta$ and as $Q_i$ is a monomial, it should be of the form
$m_1^km_2^{\delta-k}$ for some $k\in \{0,\ldots,\delta\}$.
Conversely, any set consisting of such monomials is clearly a
reducibility monomial site of $P$.

Assume next that $\psi$ is not a monomial.
From the preliminary discussion $\ell = 1$ and $Q_1=m^d$.  In
particular, $P$ and $m$ are relatively prime.  It follows from lemma
\ref{lem:monunicity}  (\ref{it:monuni2}) that $P=\psi^{d''}$ with $\psi$ homogeneous in
$m_1$ and $m_2$ and $d''\geq 2$. In particular this can only occur if $P$ is a 
pure power. Thus we are done with case (1) of theorem \ref{th:mon} (general)
where $P$ being a pure power is excluded.  If $P$ is a pure power, what  we have 
obtained is contained in  possibilities (1) and (2) 
from theorem \ref{th:mon} (addendum 2). 

\penalty -3000
\subsubsection{2nd case:} $P$ is not homogeneous in two monomials.

In this case $\psi$ is not a monomial and the desired conclusions --- that is, on one hand, case (2)  of theorem \ref{th:mon} (general) and on the other hand that only possibility (2) can occur
apart from possibilities (1) and (3) in theorem \ref{th:mon} (addendum 2) --- are part
of the preliminary discussion.


\section{Specialization}\label{sect:specialization}
In this section we explain how irreducibility properties of $F(\xbar,\lambdabar)$ can be preserved 
by specialization of  the variables $\lambda_i$ in $K$. This is the last stage
towards the results stated in the introduction. 

\subsection{Using Stein like results}

\begin{proposition} \label{prop:specialisation}
  Assume $F(\xbar,\lambdabar)=P+\lambda_1Q_1+ \cdots +
  \lambda_\ell Q_\ell$ is irreducible in $\Klambdabar [\xbar]$ (that is $\{Q_1,\ldots,
  Q_\ell\}$ is not a reducibility monomial site of $P$). Then
  for every $i=1,\ldots,\ell$, the set of $\lambda\point_i \in K$
  such that $P+\lambda_1Q_1+ \cdots + \lambda_{i-1}Q_{i-1}+
  \lambda\point_iQ_i + \lambda_{i+1}Q_{i+1}+ \cdots +\lambda_\ell
  Q_\ell$ is reducible in
  $\overline{K(\lambda_1,\ldots,\lambda_{i-1},\lambda_{i+1},\ldots,
    \lambda_\ell)}[\xbar]$ is finite and of cardinality $< \deg(P)^2$.
    \vskip 1mm
    
Consequently, for every $\lambda_1\point \in K$ but in a finite set of cardinality $<\deg(P)^2$, for every $\lambda_2\point \in K$ but in a finite set of cardinality $<\deg(P)^2$ (depending on $\lambda_1\point$),..., for every $\lambda_\ell\point \in K$ but in a finite set of cardinality $<\deg(P)^2$ (depending on $\lambda_1\point,\ldots,\lambda\point_{\ell-1}$), the polynomial $P+\lambda\point_1 Q_1+\cdots + \lambda\point_\ell Q_\ell$ is irreducible in $K[\xbar]$.

\end{proposition}

\begin{remark}
  The assumption ``$P+\lambda_1Q_1+ \cdots + \lambda_\ell Q_\ell$ 
  irreducible in $\Klambdabar [\xbar]$'' holds if it holds for a smaller
  $\ell$, in particular if $P$ itself is irreducible in $K[\xbar]$.
  This follows immediately from the equivalence of (1) and (3)
   in the Bertini-Krull theorem.
\end{remark}

\begin{proof}[Proof of proposition \ref{prop:specialisation}]
  With no loss of generality we may assume $i=1$ in the first part. Set
  $G=P+\lambda_2Q_2+ \cdots + \lambda_\ell Q_\ell$ and
  $L=\overline{K(\lambda_2,\ldots,\lambda_\ell)}$. By hypothesis,
  $G+\lambda_1 Q_1$ is irreducible in
  $\overline{L(\lambda_1)}[\xbar]$.  From the generalization of
  Stein's theorem to general pencils of hypersurfaces $P+\lambda Q$
(and not just the curves $P+\lambda$)
given in \cite{Bo} (relying on \cite{Ru}, \cite{Lo} and \cite{Na}), the set of $\lambda\point
  \in F$ such that $G+\lambda\point Q_1$ is reducible in $L[\xbar]$
  is finite and of cardinality $< \deg(P)^2$. The second part is an easy induction.
\end{proof}

\subsection{Proof of the results from the introduction}

\subsubsection{Proof of theorem \ref{thm:typical}} 
Due to the assumptions on the monomials of $P$ and $Q$, $Q$ 
cannot be a reducibility monomial in the homogeneous case (1) 
from theorem \ref{th:mon} (general) nor in possibility (1) from
theorem \ref{th:mon} (addendum 2). The monomial $Q$
not being a pure power forbids condition (2) from theorem
\ref{th:mon} (addendum 1) (with $\ell=1$ and $Q_1=Q$) to happen
and $Q$ to be a reducibility monomial in case (2) from
theorem  \ref{th:mon} (general) 
and in possibilities 
(2) and (3) from theorem \ref{th:mon} (addendum 2).  
Therefore $P+\lambda Q$ is irreducible
in $\Klambdar[\xbar]$. Apply then proposition
\ref{prop:specialisation} to complete the proof of theorem
\ref{thm:typical}.

\subsubsection{Proof of theorem \ref{thm:typical3}}
Assume as in theorem \ref{thm:typical3} that $P$ is not of the form
$h(m,\psi)$ with $h\in K [u,v]$ homogeneous of degree $\geq 2$,
$\psi\in K[\xbar]$ and $m$ a monomial dividing $Q$.  In particular
$P$ is not a pure power (for otherwise $P$ is of this form with 
$h(u,v)=v^d$ for some $d>1$ and $m=1$). We show below that assuming $Q$
is a reducibility monomial of $P$ leads to a contradiction.

The homogeneous case (1) from theorem \ref{th:mon} (general) can be 
ruled out as follows. If this case occured, then by assumption neither $m_1$ 
nor $m_2$ could divide $Q$ but this is not possible in view of the form of the
reducibility monomial sites in this case. 

The case $P$ is a monomial can also be excluded: condition (2) from 
 theorem \ref{th:mon} (addendum 1) (with $\ell=1$ and $Q_1=Q$) cannot hold
since $P$ is not a pure power. 
 
The remaining possibility (2) from theorem \ref{th:mon} (general) cannot 
happen either since in this case $P$ should be of the form
$h(m,\psi)$ as above and $Q=m^d$ (and so $m$ divides $Q$).

Conclude $Q$ is not a reducibility monomial of $P$, that is, $P+\lambda Q$
 is irreducible in $\Klambdar[\xbar]$, and apply proposition
\ref{prop:specialisation} to complete the proof of theorem
\ref{thm:typical3}.

\subsubsection{Proof of theorem \ref{thm:typical2}}  \label{proof theorem 1.3}
Here $\ell \geq 2$. The
reducibility monomial sites of cardinality $\ell$ can only occur in the homogeneous
cases from theorem \ref{th:mon} or in characteristic $p>0$. But these 
possibilities are ruled out by the assumptions. Therefore
$P+\lambda_1 Q_1+\cdots +  \lambda_\ell Q_\ell$ is irreducible in
$\Klambdabar[\xbar]$. Apply then the classical Bertini-Noether theorem
\cite[proposition 8.8]{FrJa} or alternatively proposition \ref{prop:specialisation} 
to conclude the proof.

\subsection{Further consequences}

We give below some variations around Stein's theorem which 
can be deduced from our results.

\begin{corollary} \label{cor:funny}
Let $P \in K[x_1,\ldots,x_n]$ be a polynomial in $n\geq 2$ variables and with
  coefficients in the algebraically closed field $K$.

\begin{enumerate}
\item  If $P$ is not a composed polynomial then
  $P(x_1,\ldots,x_n)+\lambda\point$ is irreducible for all but at most $\deg(P)-1$ values of
  $ \lambda\point \in K$.
  \vskip 1mm
   
\item If $P\notin K[x_1]$ and is not  divisible by $x_1$,
   then $P(x_1,\ldots,x_n)+\lambda\point x_1$
  is irreducible for all but at most $\deg(P)^2-1$ values of
  $ \lambda\point \in K$.
   \vskip 1mm
  
\item If $P\in K[x_2,\ldots,x_n]$ is not a pure power and $e$ is
an integer such that $0<e\le \deg P$ then $P(x_2,\ldots,x_n)+\lambda\point x_1^e$ 
is irreducible for all but at most $\deg(P)^2-1$ values of $ \lambda\point \in K$.
 \vskip 1mm

\item If $n=2$ and $P(x,y) \in K[x,y]$ is homogeneous of degree $d>1$ but is not a pure power and $Q=x^i y^j$
is a monomial of degree $i+j<d$ and relatively prime to $P$, then $P(x,y)+\lambda\point x^i y^j$ is irreducible for all but at most $\deg(P)^2-1$ values of $ \lambda\point \in K$.

\end{enumerate}
\end{corollary}

\begin{proof}

(1) This is the special case $Q=1$ of theorem \ref{thm:typical3} (see the comment
after theorem \ref{thm:typical3}). The bound for the number of exceptional values $\lambda\point$ is obtained by using Stein's theorem \cite{St} instead of the general bound from [Bo] as in proposition
\ref{prop:specialisation}.

(2) Suppose that $P(x_1,\ldots,x_n)+\lambda x_1$ is reducible in $\Klambdar[\xbar]$. As $x_1$ is not a pure power, it follows from theorem \ref{th:mon} that $P=h(m_1,m_2)$ for some homogeneous polynomial $h\in K[u,v]$ of degree $d>1$ and some monomials $m_1$ and $m_2$ and that $x_1 = m_1^km_2^{d-k}$ for some $k\in \{0,\ldots,d\}$. Then we have necessarily $\{m_1,m_2\} = \{1,x_1\}$. But then $P=h(m_1,m_2)$ contradicts the assumption $P\notin K[x_1]$. Thus $P(x_1,\ldots,x_n)+\lambda x_1$ is irreducible in $\Klambdar[\xbar]$ and the result follows from proposition \ref{prop:specialisation}.

(3)  We first show that $P(x_2,\ldots,x_n)+\lambda x_1^e$ is irreducible in $\Klambdar[\xbar]$. 
From theorem \ref{th:mon} we need to exclude the two following situations.

\begin{enumerate}
 \item $P=h(m_1,m_2)$ for some homogeneous polynomial $h\in K[u,v]$ of degree $d>1$ and some relatively prime monomials $m_1$ and $m_2$ and $x_1^e = m_1^k m_2^{d-k}$ with $0\le k \le d$. If $0<k<d$ then necessarily one of the two monomials, say $m_1$, is constant and $m_2$ is a pure power of $x_1$. But then $P=h(m_1,m_2)$ contradicts the assumption $\deg_{x_1}(P)=0$.  If $k=0$, $m_2$ is a pure power of $x_1$ but then $P=h(m_1,m_2)$ is possible only if $P=m_1^d$ (for otherwise $\deg_{x_1}(P)>0$), which is excluded as $P$ is not a pure power. The case $k=d$ is similar.
  
\item  $P=h(m,\psi)$ for some homogeneous polynomial $h\in K[u,v]$ of degree $d>1$ and $x_1^e = m^d$. Then $m$ is a pure power of $x_1$ and as above $P=h(m,\psi)$ is possible only if $P=\psi^d$ (for otherwise $\deg_{x_1}(P)>0$), which is excluded as $P$ is a not pure power.

  \end{enumerate}

\noindent  The result  follows then from proposition \ref{prop:specialisation}.
\vskip 1mm
  
(4) Irreducibility of $P(x,y)+\lambda x^i y^j$ in $\Klambdar[x,y]$ readily follows from theorem \ref{th:mon} (general \& addendum 1):  just note $P$ is homogeneous in the two monomials $m_1=x$ and $m_2=y$, which are relatively prime, of degree $<\deg(P)$ and such that $m_1+\lambda m_2$ is irreducible in $\Klambdar[x,y]$. Apply then proposition \ref{prop:specialisation} to complete the proof.
\end{proof}


\end{document}